\documentclass[a4paper,12pt]{article} 
\usepackage{amsmath, amsthm, amssymb}
\usepackage{url}

\usepackage{physics}
\usepackage{graphicx,lipsum}
\graphicspath{ {./downloads/} }

\usepackage[colorlinks,citecolor=red,urlcolor=blue,bookmarks=false,hypertexnames=true]{hyperref}
\usepackage{tikz}
\usetikzlibrary{calc}
\usetikzlibrary{shapes}
\usepackage[autostyle]{csquotes}
\makeatletter

\usepackage[colorlinks]{hyperref}
\usepackage[nameinlink,capitalize]{cleveref}
\newtheorem{theorem}{Theorem}[section]
\newtheorem{corollary}[theorem]{Corollary}
\newtheorem{lemma}[theorem]{Lemma}
 
 \newtheorem{remark}[theorem]{Remark}
\newtheoremstyle{named}{}{}{\itshape}{}{\bfseries}{.}{.5em}{\thmnote{#3's }#1}
\theoremstyle{named}

\theoremstyle{definition}
\newtheorem{definition}[theorem]{Definition}
\newtheorem{example}[theorem]{Example}

\usepackage{xspace}
\usepackage[margin=1.0in]{geometry}

\usepackage{titlesec}

\usepackage{mathtools}

\DeclarePairedDelimiterX{\inp}[2]{\langle}{\rangle}{#1, #2}
\titleformat{\chapter}
  {\Large\bfseries} 
  {}                
  {0pt}            
  {\huge}

\begin{document}

\begin{center}
\fontsize{13pt}{10pt}\selectfont
    \textsc{\textbf{GRADED PSEUDO WEAKLY PRIME SPECTRUM OF GRADED TOPOLOGICAL MODULES}}
    \end{center}
\vspace{0.1cm}
\begin{center}
   \fontsize{12pt}{10pt}\selectfont
    \textsc{{\footnotesize Tamem Al-shorma$n$, Malik Bataine$h$, Melis Bolat$^\ast$, and Bayram Ali Ersoy }}
\end{center}
\vspace{0.2cm}

\begin{abstract}
  		In this study, we introduce graded pseudo weakly prime submodules of G-graded R-modules, which are an extension of graded weakly prime ideals over G-graded rings. On the graded spectrum of graded pseudo weakly prime submodules, we investigate the Zariski topology. Different aspects of this topological space are investigated, and they are linked to the algebraic properties of the G-graded R-modules under study.
\end{abstract}

\section{introduction}
In this paper, we introduce a Zariski topology on the graded pseudo weakly prime spectrum $\Upsilon_M$ of graded pseudo weakly prime submodules of a certain $G$-graded $R$-module $M$ over $G$-graded ring $R$, and investigate the relation between the properties of $M$ and the topological space we obtain. The relation between the Zariski topology established on the graded prime spectrum of a $G$-graded ring $R$ and the ring theoretic features of $R$ in \cite{ozkiirisci2013graded} is the motivation for this study.
	\\
	
	Let G be an abelian group with identity e. A ring R is called a G-graded ring if $ R= \bigoplus\limits_{g \in G} R_g$   with the property $R_gR_h\subseteq R_{gh}$ for all $g,h \in G$, where $R_g$ is an additive subgroup of R for all $g\in G$. The elements of $R_g$ are called homogeneous of degree g. If $x\in R$, then $x$ can be written uniquely as $\sum\limits_{g\in G} x_g$, where $x_g$ is the component of $x$ in $R_g$. The set of all homogeneous elements of R is $h(R)= \bigcup\limits_{g\in G} R_g$. Let P be an ideal of a G-graded ring R. Then P is called a graded ideal if $P=\bigoplus\limits_{g\in G}P_g$, i.e, for $x\in P$ and  $x=\sum\limits_{g\in G} x_g$ where $x_g \in P_g$ for all $g\in G$. An ideal of a G-graded ring is not necessarily a graded ideal (see \cite{abu2019graded}).  
	\\
	
	Let P be a proper graded ideal of R. The graded radical of P is denoted by $Grad(P)$ and defined as follows:
	\begin{center}
		{\small $Grad(P)=\Big{\{} x= \sum\limits_{g\in G} x_g \in R$ : for all $g\in G$, there exists a $n_g\in \mathbb{N}$ such that ${x_g}^{n_g} \in P \Big{\}}$}.  
	\end{center}
	Note that $Grad(P)$ is always a graded ideal of R.
	\\
	
	Let R be a G-graded ring. A left R-module M is said to be a graded R-module if there exists  a family of additive subgroups $\{M_g\}_{g \in G}$ of M such that $M= \bigoplus\limits_{g \in G} M_g$ with the property $M_gM_h\subseteq M_{gh}$ for all $g,h \in G$. The set of all homogeneous elements of M is $h(M)= \bigcup\limits_{g\in G} M_g$. Note that $M_g$ is an $R_e$-module for every $g\in G$. A submodule P of M is said to be graded submodule of M if $P=\bigoplus\limits_{g\in G}P_g$. Here, $P_g$ is called g-component of P (see \cite{gordon1982graded}).
	\\
	
	Graded weakly prime ideals have been introduced by S. E. Atani in \cite{atani2006grade}; a proper graded ideal $I$ of a $G$-graded ring $R$ is said to be graded weakly prime ideal of $R$ if whenever $a,b\in h(R)$ such that $0\neq ab\in I$, then either $a\in I$ or $b\in I$.
	\\
	
	In \cite{atani2006graded}, S. E. Atani defined graded weakly prime submodules; a proper graded submodule $P$ of a $G$-graded $R$-module $M$ is said to be graded weakly prime submodule of $M$ if whenever $r\in h(R)$ and $m\in h(M)$ such that $0\neq rm\in P$, then either $m\in P$ or $r\in (P:_{R}M)$ where $(P:_{R}M)$ =$\{r\in R:rM\subseteq P\}$ is a graded ideal of $R$, (see \cite[Lemma 2.1]{farzalipour2012union}).
	\\
	
	Let $R$ be a $G$-graded ring, $M$ a graded $R$-module and $P$ a graded submodule of $M$. $(P:_{R}M)$ is defined as $(P:_{R}M)=\{r\in R:rM\subseteq P\}.$ It is shown in \cite{atani2006grad} that if $P$ is a graded submodule of $M$, then $(P:_{R}M)$ is a graded ideal of $R$. The annihilator of $M$ is defined as $(0:_{R}M)$ and is denoted by $Ann_{R}(M)$. A graded submodule $P$ of $M$ is said to be a graded maximal if $P\neq M$ and if there is a graded submodule $L$ of $M$ such that $P\subseteq L\subseteq M$, then $L=P$ or $L=M$ \ (see \cite{nastasescu2004methods}).
	\\
	
	In this topic, there are numerous studies in the literature. In \cite{tekir2009zariski}, U. Tekir studied the Zariski topology on the prime spectrum of a module over non-commutative rings. Graded 2-absorbing quasi primary ideal was introduced in \cite{uregen2019graded}  as a generealization of graded prime ideals. Another generealization of graded prime ideals, graded 1‑absorbing prime ideal, was defined in \cite{abu2021graded}. In \cite{yildiz2020s}, S-Zariski topology was constructed. 
	\\
	
	In this paper, we introduce the concept of graded pseudo weakly prime submodules of $G$-graded $R$-modules which is a generalization of the graded weakly prime ideals over $G$-graded rings. Let $R$ be a $G$-graded ring, $M$ a $G$-graded $R$-module and $P$ a graded submodule of $M$. A graded submodule $P$ of a $G$-graded $R$-module $M$ is a graded pseudo weakly prime if $(N:_{R}M)$ is a graded weakly prime ideal of $R$.
	\\
	
	We construct the Zariski topology on the graded spectrum of graded pseudo weakly prime submodules. We clarify the relationship  between the properties of this topological space and the algebraic properties of the $G$-graded $R$-modules.

	\section{Graded Pseudo Weakly Prime Submodules}
	In this section, we investigate graded pseudo weakly prime submodules, with a variety of outcomes.
	
	\begin{definition}\label{def1}
		Let M be a G-graded R-module and P be a graded submodule of M.
		\\
		
		(1) P is a graded pseudo weakly prime submodule of M if $(P:_{R}M)$ is a graded weakly prime ideal of $R$.
		\\
		
		(2) The graded pseudo weakly prime spectrum of a G-graded $R$-module $M$ is defined to be the set of all graded pseudo weakly prime submodules of $M$, and denoted by $\Upsilon_{M}$. For any graded weakly prime ideal $I$ of $R$, the collection of all graded pseudo weakly prime submodule $P$ of $M$ with $(P:_{R}M)=I$ is denoted by $\Upsilon_{M,I}$.
		\\
		
		(3) A graded submodule P of M is defined $\chi^M(P)=\{I\in \Upsilon_{M}:P\subseteq I\}$. We use $\chi(P)$ instead of $\chi^M(P)$ for nonambiguity.
		\\
		
		(4) Suppose that the set of all graded pseudo weakly prime submodules of G-graded R-module $M$ is not empty set (i.e $\Upsilon_{M}\neq \emptyset $). The map $\varphi :\Upsilon_{M}\rightarrow GWSpec(R/Ann(M))$ where $GWSpec(R)$ denotes the set of all graded weakly prime ideals of $R$, defined by $\varphi(P) = (P:_{R}M)/Ann(M)$ is called the natural map of $\Upsilon_{M}$.
		\\
		
		(5) A G-graded $R$-module $M$ is said to be G-graded pseudo weakly primeful R-module if either $M=\{0\}$ or $M\neq \{0\}$ and $\varphi $ is surjective. If  $\varphi $ is injective, then $M$ is called G-graded pseudo weakly injective R-module.
	\end{definition}
	
	The graded weakly prime ideals of the $G$-graded ring $R$ and the graded pseudo weakly prime submodules of the $G$-graded $R$-module $M$ are the same, according to our definition. This demonstrates that the graded pseudo weakly prime submodule is a $G$-graded $R$-module version of the graded weakly prime ideal.
	\\
	
	Recall that, a graded submodule $P$ of a $G$-graded $R$-module $M$ is said to be graded prime submodule if $P \neq M$ and whenever $x \in h(R)$ and $m \in h(M)$ with $0\neq xm \in P$, then either $m\in P$ or $y \in (P:_RM)$ ( see \cite{atani2006graded}). Every graded weakly prime submodule of G-graded R-module M is graded pseudo weakly prime submodule of M, since $(P:_RM) \in GWSpes(R)$. However, the converse is not true in general.
	
	\begin{example}\label{ex1}
		Let $M=\mathrm{Z}\oplus \mathrm{Z}$ be a G-graded $\mathrm{Z}$-module, where $M_0=M$ and $M_g=\{0\}$ for all $g\in G$ and $P = (2,0)\mathrm{Z}$ is a graded submodule of M generated by $(2,0) \in h(M)$. Thus $(P:_RM) = \{0\} \in GWSpec(\mathrm{Z})$ (i.e. $P \in \Upsilon_M$). Despite the fact that P is not a weakly prime submodule of M. Therefore, a graded pseudo weakly prime submodule need not be graded weakly prime submodule (i.e. $GWSpec(M) \subset \Upsilon_M$), where GWSpec(M) is the set of all graded weakly prime ideals of M.    
	\end{example}
	
	Note that, a $G$-graded $R$-module $M$ is called $G$-graded multiplication $R$-module if for every graded submodule $P$ of $M$, $P=MI$ for some graded ideal $I$ of $R$ (see \cite{escoriza1998multiplication}). Clearly, every $G$-graded multiplication $R$-module is $G$-graded pseudo weakly injective $R$-module.
	\\
	
	To prove that let $P$, $I$ be two graded pseudo weakly prime submodules of $M$ such that $\varphi (P)=\varphi (I)$, so $(P:_{R}M)+(0:_{R}M)=(I:_{R}M)+(0:_{R}M)$ then $(P:_{R}M)-(I:_{R}M)\subseteq (0:_{R}M)$ (i.e. $(P:_{R}M)M-(I:_{R}M)M=\{0\}$), since $M$\ is $G$-graded
	multiplication $R$-module then $P-I=\{0\}$, we have $P=I$, so $M$ is $G$-graded pseudo weakly injective $R$-module.
	\\
	
	We will show that every $G$-graded pseudo weakly injective $R$-module is $G$-graded multiplication $R$-module provided that the module is graded finitely generated.
	
	\begin{corollary}\label{coro1}
		Let $M$  be a G-graded finitely generated $R$-module. If $M$ is G-graded pseudo weakly injective R-module, then $M$ is G-graded multiplication R-module.
		\\
		\\
		\textbf{Proof.} Let $M$ be G-graded pseudo weakly injective R-module. If $|\Upsilon_{M,I}|\leq 1$ for every graded maximal ideal $I$ of $R$ and then by Lemma 2.9 in \cite{abu2021zariski}, $M/IM$ \ is a graded simple and so graded cyclic and then by Corollary 2.8 in \cite{abu2021zariski}, $M$ is a G-graded multiplication $R$-module.
	\end{corollary}
	
	\begin{remark}\label{rem1}
		Let $P\subseteq \Upsilon_{M}$. The intersection of all elements in $P$ is denoted by ${\eta}(P)$.
	\end{remark}
	
	\begin{definition}\label{def2}
		Let M be a G-graded R-module and P be a graded submodule of M. Then the following statements are true:
		\\
		
		(1) P is said to be a graded pseudo weakly semiprime submodule of M if $P$ is an intersection of graded pseudo weakly prime submodules of $M$.
		\\
		
		(2) If P is a graded pseudo weakly prime submodule of M. Then $P$ is said to be a graded weakly extraordinary submodule of M if whenever $I$ and $J$ are graded pseudo weakly semiprime submodules of $M$ such that $\{0\}\neq I\cap J\subseteq P$, then either $I\subseteq P$ or $J\subseteq P$.
		\\
		
		(3) M is said to be a G-graded weakly topological R-module if $\Upsilon_{M}=\emptyset $ or every graded pseudo weakly prime submodule of $M$ is a graded weakly extraordinary submodule.
		\\
		
		(4) A graded pseudo weakly prime radical of $P$ will be denoted by $GPWrad(P)$ and is defined by the intersection of all graded pseudo weakly prime submodules of $M$ containing $P$ (i.e. $GPWrad(P)={\eta}(\chi(P))=\bigcap\limits_{I\in \chi(P)}I$).
		\\
		
		(5) If $\chi(P)=\emptyset $, then we say that $GPWrad(P)=M$. Also, $P$ is said to be a graded pseudo weakly prime radical if $GPWrad(P)=P$.
	\end{definition}

	\begin{lemma}\label{lem1}
		Let $M$ be a G-graded $R$-module. Then $M$ is a G-graded weakly topological R-module if and only if $\chi(P)\cup \chi(I)=\chi(P\cap I)$ for every graded pseudo weakly semiprime submodules $P$ and $I$ of $M$.
		\\
		\\
		\textbf{Proof.} Suppose that $M$ is a G-graded weakly topological $R$-module. Let $P$ and $I$ be graded pseudo weakly semiprime submodules of $M$. Clearly, $\chi(P)\cup \chi(I)\subseteq \chi(P\cap I)$. To prove that let $J\in \chi(P)\cup \chi(I)$ such that $J$ be a graded pseudo weakly prime submodules of $M$ then $J\in \chi(P)$ or $J\in \chi(I)$, so $P\subseteq J$ or $I\subseteq J$. Hence $P\cap I\subseteq J$, also $J \in \chi(P\cap I)$, so we are done. Let $K\in \chi(P\cap I)$. Then $P\cap I\subseteq K$ and then by assumption (since $M$ is a G-graded weakly topological R-module), either $P\subseteq K$ or $I\subseteq K$ (i.e. either $K\in \chi(P)$ or $K\in \chi(I)$)\ then $K\in \chi(P)\cup \chi(I)$ which proves that $\chi(P\cap I)\subseteq \chi(P)\cup \chi(I)$. Hence, $\chi(P)\cup \chi(I)=\chi(P\cap I)$. Conversely, let $J$ be a graded pseudo weakly prime submodule of $M$. Assume that $I$ and $L$ are graded pseudo weakly semiprime submodules of $M$ such that $P\cap L\subseteq J$. Then $J\in \chi(P\cap I)=\chi(P)\cup \chi(I)$ by assumption. So, $J\in \chi(P)$ or $J\in \chi(I)$ that is either $P\subseteq J$ or $I\subseteq J$. Hence, $M$ is a G-graded weakly topological $R$-module.
	\end{lemma}
	
	\begin{theorem}\label{thm1}
		Let $R$ be a $G$-graded ring and let $P$ be a graded ideal of $R$. Then the following statements are equivalent:
		
		$(i)$ $P$ is a graded weakly prime ideal of $R$.
		
		$(ii)$ For each graded ideals $I_{1}$, $I_{2}$ of $R$ with $\{0\}\neq I_{1}I_{2}\subseteq P$, we have either $I_{1}\subseteq P$ or $I_{2}\subseteq I$.
		\\
		\\
		\textbf{Proof.} $(i)\Rightarrow (ii)$\ Let $P$  be a graded
		weakly prime ideal and let $I_{1}$, $I_{2}$ be two graded ideals of $R$ with $ \{0\}\neq I_{1}I_{2}\subseteq P$. Suppose that $I_{1}\nsubseteq P$ and $I_{2}\nsubseteq P$. Then there exist $x\in(h(R)\cap I_{1})-P$ and $y\in (h(R)\cap I_{2})-P$. Hence  $x\in P$ or $y \in P$ which is a contradiction. Therefore, $I_{1}\subseteq P$ or $I_{2}\subseteq I$.
		\\
		$(ii)\Rightarrow (i)$ Let $x,y \in h(R)$ with $0\neq xy\in P$.
		Let $I_{1}=Rx$ and $I_{2}=Ry$. Then $I_{1}$, $I_{2}$ are graded ideals of $R$ and $\{0\}\neq \ I_{1}I_{2}\subseteq P$. Hence $I_{1}\subseteq P$ or $I_{2}\subseteq P$ (i.e. $x\in P$ or $y\in P)$.
		
	\end{theorem}

	\begin{theorem}\label{thm2}
		Let M be a G-graded R-module. If M is a G-graded multiplication R-module then M is a G-graded weakly topological R-module.
		\\
		\\
		\textbf{Proof.} Suppose that $M$ is a G-graded multiplication $R$-module. Let $N$ be a graded pseudo weakly prime submodule of $M$. Assume that $K$ and $L$ are two graded pseudo weakly semiprime submodules of $M$ such that $\{0\}\neq K\cap L\subseteq N$. Then there exist graded ideals $I$ and $J$ of $R$ such that $\chi(K)=\chi(IM)$ and $\chi(L)=\chi(JM)$. Suppose that $K=\cap_{i\in \Delta }P_{i}$ for some graded pseudo weakly prime submodules $\{P_{i}\}_{i\in \Delta }$. Now, for every $i\in \Delta $, $P_{i}$ $\in
		\chi(K)\subseteq \chi(K)\cup \chi(L)$ $=\chi(IM)\cup \chi(JM)\subseteq \chi((I\cap J)M)$ by Lemma \ref{lem1}. This implies that $(I\cap J)M\subseteq \cap _{i\in \Delta}P_{i}=K$. Similarly, $(I\cap J)M\subseteq L$. So, $(I\cap J)M\subseteq K$ $\cap $ $L\subseteq N$ (i.e. $I\cap J\subseteq (N:_{R}M))$, also since $(N:_{R}M)$ is a graded weakly prime ideal of $R$ and $IJ\subseteq I\cap J$, so $\{0\}\neq \ IJ$\ $\subseteq (N:_{R}M)$\ then either $I\subseteq
		(N:_{R}M) $ or $\ J\subseteq (N:_{R}M)$, which implies that $K=IM\subseteq N$ or $L=JM\subseteq N$. Hence, $M$ is a G-graded weakly topological $R$-module.
	\end{theorem}
	
	\begin{lemma}\label{lem2}
		Let M be a G-graded R-module and P be a graded submodule of M. If $M$ is a G-graded weakly topological $R$-module, then $M/P$ is a G-graded weakly topological $R$-module.
		\\
		\\
		\textbf{Proof.} Any graded pseudo weakly semiprime submodule of $M/P$ has the form $I/P$ where $I$ is a graded pseudo weakly semiprime submodule of $M$ containing $P$. Assume that $I_{1}/P$ and $I_{2}/P$ are two graded pseudo weakly semiprime submodules of $M/P$, also since $M$ is a G-graded weakly topological R-module we can apply Lemma \ref{lem1}, so for any graded pseudo weakly semiprime submodules $I_{1}$ and $I_{2}$ of $M$  we have $\chi(I_{1})\cup \chi(I_{2})=\chi(I_{1}\cap I_{2})$, also $I_{1}$ $\subseteq P$ \ and $I_{2}\subseteq P$ \ then $I_{1}/P\ =I_{1}$ and $I_{2}/P=$ $I_{2}$, so $\chi(I_{1}/P)\cup \chi(I_{2}/P)=\chi(I_{1}/P\cap I_{2}/P).$ Hence $M/P$  is a G-graded weakly topological R-module.
	\end{lemma}

	\section{Graded Pseudo Weakly Prime Spectrum Of Graded Topological
		Modules}
	
	Now, we introduce and study a topological on the graded pseudo weakly prime spectrum and give several results concerning it as a generalization of a topological on graded pseudo prime spectrum in \cite{abu2021zariski}. This topology is called Zariski topology.
	\\
	
	Also, in this section, we assume that $M$ is always $G$-graded weakly topological $R$-module. Then for each graded submodule $P$ of $M$ we define the variety of $P$ \ by $\chi(P)=\{I\in \Upsilon_{M},P\subseteq I\}$, so this variety satisfies
	the topology axioms for closed set, since $\phi =\chi(M)$, $\Upsilon_{M}=\chi(\{0\})$ and for any family of graded submodules $\{Pi\}_{i\in \Delta }$ of \ $M$,  $\bigcap\limits _{i\in \Delta }\chi(Pi)=\chi(\sum\limits_{i\in \Delta }Pi)$. Also, for any
	graded submodules $P$ and $J$ of $ M$, $\chi(P)\cup \chi(J)=\chi(P\cap J)$. Thus, if $\xi (M)$ denotes the collection of all subsets $\chi(P)$ of $\Upsilon_{M}$ \ (i.e. $\xi (M)=\{\chi(P)$, $P$ is a graded submodule of $M$\} contain $\Upsilon_{M}$ \ and empty set) closed under arbitrary intersections and it is closed under finite union, then $\xi (M)$ satisfies the axioms of a topological space for the closed subsets. This topology is called the Zariski topology.
	\\
	
	The topological aspects of this topology are investigated, and some results about the relationship between algebraic properties of graded topological modules and topological properties of the Zariski topology on the graded pseudo weakly prime spectrum of them are presented.
	\\
	
	We define $\chi^{R}(P)=\{I\in \Upsilon_{R}:P\subseteq I\}$. Moreover, for every graded ideal $J\in \chi^{R}(Ann(M))$, $R^{\star }$ and $K^{\star }$ will denote $R/Ann(M)$ and $K/Ann(M),$ respectively.
	
	\begin{theorem}\label{thm3}
		Let $M$ be a G-graded pseudo weakly primeful $R$-module such that $\Upsilon_{M}$ is connected. Then $\Upsilon_{R^{\star }}$ is connected.
		\\
		\\
		\textbf{Proof.} Since the natural map $\varphi:\Upsilon_{M}\rightarrow GWSpec(R/Ann(M))$ such that $\varphi (P)=(P:_{R}M)/Ann(M)$ is surjective. It is sufficient to demonstrate that $\varphi$ is a continuous map in terms of the Zariski topology. To prove it, let $K$ be a graded ideal of $R$ containing $Ann(M)$ and let $P\in \varphi ^{-1}(\chi^{R^{\star }}(K^{\star }))$. Then there exists $I^{\star}\in \chi^{R^{\star }}(J^{\star }))$ such that $\varphi (P)=I^{\star }$ and then $I=(P:_{R}M)\supseteq J$ which implies that $KM\subseteq P$ and hence $P\in\chi(KM)$. Let $L\in \chi(KM)$. Then$(L:_{R}M)\supseteq (KM:_{R}M)\supseteq J$ and then $L\in \varphi ^{-1}(\chi^{R^{\star }}(K^{\star }))$. So, $\varphi^{-1}(\chi^{R^{\star }}(K^{\star }))=\chi(KM)$ which means that $\varphi $ is continuous.
	\end{theorem}
	
	\begin{remark}\label{rem2}
		Let$\ M$ be a G-graded $R$-module and $W$ be any subset of $\Upsilon_{M}$, we will denote the closure of $W$ in $\Upsilon_{M}$ for the zariski topology by $Cl(W)$.
	\end{remark}
	
	We show in the next theorem that if the topological space $\Upsilon_M$ is a $T_0$-space, we can derive some properties for the graded pseudo weakly prime submodules of $M$.
	
	\begin{theorem}\label{thm4}
		Let M be a G-graded R-module, $P \subseteq \Upsilon_M$ and $J \in \Upsilon_{M,I}$ for some $I \in GWSpec(R)$.
		\\
		
		(1) $Cl(P)=\chi({\eta}(P))$ and hence $P$ is closed if and only if $P=\chi({\eta}(P))$. In particular, $Cl(\{J\})=\chi(J)$.
		\\
		
		(2) If $\{0\}\in P$, then $P$ is dense in $\Upsilon_{M}$.
		\\
		
		(3) $\Upsilon_{M}$ is a $T_{0}$- space.
		\\
		\\
		\textbf{Proof.} (1) Clearly, $P\subseteq \chi({\eta}(P))$ since ${\eta}(P)\subseteq P$. Let $\chi(N)$ be any closed subset of $\Upsilon_{M}$ containing $P$. Since ${\eta}(P)\supseteq{\eta}(\chi(P)),$ then $\chi({\eta}(P))\subseteq \chi({\eta}(\chi(N)))=\chi(GPWrad(N))=\chi(N)$ which means that $\chi({\eta}(P))$ is the smallest closed subset of $\Upsilon_{M}$ containing $P$. Hence, $Cl(P)=\chi({\eta}(P))$.
		\\
		
		(2) It is clear by (1), $Cl(P)=\chi({\eta}(P))$, since $\{0\}\in P$ then ${\eta}(P)=$ $\{0\}$ so $Cl(P)=\chi(\{0\})=\Upsilon_{M}$ (a topology space) then $P$ is dense in $\Upsilon_{M}$.
		\\
		
		(3) The topological space is a $T_0$-space if and only if the closures of distinct points are distinct. Let $N$ and $K$ be two distinct points of $\Upsilon_{M}$. Then by (1) we have that $Cl(\{N\})=\chi(N)\neq \chi(K)=Cl(\{K\})$. Therefore, $\Upsilon_{M}$ is a $T_{0}$-space.
	\end{theorem}
	
	\begin{theorem}\label{thm5}
		Let $R$ be a $G$-graded ring, $M$ be a G-graded $R$-module, $C\in \Upsilon_{M}$ and $\Omega =\{(D:_{R}M): D\in
		\Upsilon_{M}\} \leq GWSpec(R)$. Then the set $\{C\}$ is closed in $\Upsilon_{M}$ if and only if (i) $p=(C:_{R}M)$ is a maximal element of the set $\Omega $, and (ii) $\Upsilon_{M,p}=\{C\}$, that is $\mid \Upsilon_{M,P}\mid =1$.
		\\
		\\
		\textbf{Proof.} ($\Rightarrow$) Assume that $\{C\}$ is closed in $\Upsilon_{M}$ by (1) in Theorem \ref{thm4}, we have $\{C\}=Cl(\{C\})=\chi(C)$. Let $q\in \Omega $ such that $p\subseteq q$. Then there exists $D\in \Upsilon_{M}$ such that $q=(D:_{R}M)$. Hence $(C:_{R}M)=p$ $\subseteq q=(D:_{R}M).$ Thus $D\in \chi(C)=\{C\}$ so that $D=C$ and $p=q$, (the only ideal containing p is p) then p is maximal,
		which proves $(i)$. Let $C_{0}\in \Upsilon_{M,P}$ Thus $(C_{0}:_{R}M)=p=(C:_{R}M)$ and hence $C_{0}\in \chi(C)=\{C\}$. So $C_{0}=C$. Thus, $\Upsilon_{M,P}=C$.
		\\
		
		($\Leftarrow$) We assume $(i)$ and $(ii)$, and show that $\chi(C)=\{C\}$. If $D\in \chi(C)$, then $q=(D:_{R}M)\supseteq (C:_{R}M)=p$. Therefore $(i)$ implies $q=p$ and consequently $(ii)$ implies $D=C$ \ by $(ii)$. Thus, $\chi(C)\subseteq \{C\} $. Since $C$ is a graded pseudo weakly prime submodule we have clearly $\{C\}\subseteq \chi(C)$, so that $\chi(C)=\{C\}$. By (1) in Theorem \ref{thm4}, $CL(\{C\})=\chi(C)=\{C\}.$ Therefore the set $\{C\}$ is closed in $\Upsilon_{M}$.
	\end{theorem}
	
	\begin{theorem}\label{thm6}
		Let $M$ be a G-graded $R$-module. Then $\Upsilon_{M}$ is a $T_{1}$-space if and only if every element of $\Upsilon_{M}$ is
		maximal element. (i.e. Let $R$ be a $G$-graded ring and $M$ a G-graded $R$-module, and $\Omega =\{(D:_{R}M)\mid D\in \Upsilon_{M}\}\leq GWSpec(R)$. Then $\Upsilon_{M}$ \ is a $T_{1}$
		-space if and only if $(i)$ $p=(C:_{R}M)$ is a maximal element of the set $\Omega $, $\forall C\in \Upsilon_{M}$ and $(ii)\mid \Upsilon_{M,P}\mid =1$, $\forall p\in GWSpec(R)$).
		\\
		\\
		\textbf{Proof.} If $\Upsilon_{M}$ is a $T_{1}$-space then the singleton sets are closed in $\Upsilon_{M}$, so we get $(i)$ and $(ii)$ by Theorem \ref{thm5}. Conversely, $(i)$ and $(ii)$ are equivalent so that the singleton set $\{C\}$ is closed in$\Upsilon_{M}$ for every $C\in $ $\Upsilon_{M}$, that is $\Upsilon_{M}$ is a $T_{1}$-space.
	\end{theorem}
	
	Recall that, a topological space $T$ is called irreducible if and only if for any decomposition $T = T_1 \cup T_2$ with closed subset $T_i$ of $T$ for all $i=1,2$, we have that $T_1 = T$ or $T_2=T$. If a nonempty subset $B$ of a topological space $T$ cannot be represented as the union of two proper subsets, each of which is closed in $B$, it is said to be irreducible. The empty set isn't regarded as irreducible. A maximal irreducible subset of a topological space $T$ is referred to as an irreducible component of $T$. Because every singleton subset of $\Upsilon_M$ is irreducible, so is its closure. Now, by using $(1)$ of Theorem \ref{thm4}, we obtain that.
	
	\begin{corollary}\label{coro2}
		Let $M$  be a G-graded $R$-module. Then $\chi(P)$ is an irreducible closed subset of $\Upsilon_{M}$ for every graded pseudo weakly prime submodule $P$ of $ M$.
		\\
		\\
		\textbf{Proof.} Since by (1) of Theorem \ref{thm4}, $Cl(\{P\})=\chi(P),$ also $\{P\}$ and $Cl(\{P\})$ are irreducible, then $\chi(P)$ is an irreducible closed subset of $\Upsilon_{M}$.
	\end{corollary}
	
	\begin{theorem}\label{thm7}
		Let $M$ be a G-graded $R$-module and $P\subseteq \Upsilon_{M}$. Then ${\eta}(P)$ is a graded pseudo weakly prime submodule of $M$ if and only if $P$ is an irreducible space.
		\\
		\\
		\textbf{Proof.} ($\Rightarrow$) Suppose that $W$ is irreducible. Clearly, ${\eta}(P)$ is a graded submodule of $M$. Let $I$ \ and $J$ be graded ideals of $R$ such that $0\neq IJ\subseteq ({\eta}(P):_{R}M)$. Then $0\neq IJM\subseteq {\eta}(P)$, also since $ {\eta}(P)\subseteq P$ . So, $P\subseteq \chi({\eta}(P))\subseteq \chi(IJM)=\chi(IM)\cup \chi(JM)$. Since $P$ is irreducible, either $
		P\subseteq \chi(IM)$ or $P\subseteq \chi(JM)$ and hence either ${\eta}(P)\supseteq {\eta}(\chi(IM))=GPWrad(IM)\supseteq IM$ or ${\eta}(P)\supseteq JM$ \ which implies that either $I\subseteq ({\eta}(P):_{R}M)$ or $J\subseteq ({\eta}(P):_{R}M)$, then $({\eta}(P):_{R}M)$ is graded weakly prime ideal of $R$. Hence, ${\eta}(P)$ is a graded pseudo weakly prime submodule of $M$.
		\\
		
		($\Leftarrow$) Assume that $P\subseteq P_{1}\cup P_{2}\neq 0$ for closed subsets $P_{1}$ and $P_{2}$ of $\Upsilon_{M}$. Then there exist graded submodules $N$ and $K$ of $M$ such that $P_{1}=\chi(N)$ and $P_{2}=\chi(K)$ and hence ${\eta}(P)\supseteq {\eta}(\chi(N)\cup \chi(K))={\eta}(\chi(N)\cap \chi(K))=GPWrad(N)\cap GPWrad(K)\neq 0$. Since $M$ is a G-graded weakly topological $R$-module, ${\eta}(P)$ is a graded weakly extraordinary submodule of $M$ and hence either $
		GPWrad(N)\subseteq {\eta}(P)$ or $GPWrad(K)\subseteq {\eta}(P)$. So, $P\subseteq \chi({\eta}(P))\subseteq \chi(GPWrad(N))=\chi(N)=P_{1}$ or $P\subseteq P_{2}$. Hence, $P$ is an irreducible space.
	\end{theorem}
	
	\begin{corollary}\label{coro3}
		Let M be a G-graded R-module and P be a graded submodule of M. Then the following statements are true:
		\\
		
		(1) $\chi(N)$ is an irreducible space if and only if $GPWrad(P)$ is a graded pseudo weakly prime submodule of $M$.
		\\
		
		(2) $\Upsilon_{M}$ is an irreducible space if and only if $GPWrad(\{0\})$ is a graded pseudo weakly prime submodule of $M$.
		\\
		
		(3) If $\Upsilon_{M,J}\neq \phi $  for some $J\in GWSpec(R)$, then $\Upsilon_{M,J}$ is an irreducible space.
		\\
		
		(4) If R is a quasi-local G-graded ring, then GMax(M) is an irreducible space (here the GMax(M) is the set of all graded maximal submodules of M).
		\\
		
		(5) If $\{0\}\in \Upsilon_{M}$, then $\Upsilon_{M}$ is an irreducible space.
		\\
		
		(6) If $R$ is a G-graded integral domain, then $\Upsilon_{M}$ is an irreducible space.
		\\
		\\
		\textbf{Proof.} (1) Since $GPWrad(P)={\eta}(\chi(P))$, the result is clear by Theorem \ref{thm7} (put $P=\chi(P)$, so done).
		\\
		
		(2) Take $P=\{0\}$ in (1).
		\\
		
		(3) Since $({\eta}(\Upsilon_{M,J}):_{R}M)=\bigcap\limits _{P\in \Upsilon_{M,J}}(P:_{R}M)=J\in GWSpec(R)$, so ${\eta}(\Upsilon_{M,J})$ is a graded pseudo weakly prime submodule of $M$, then$\Upsilon_{M,J}$ is an irreducible space by Theorem \ref{thm7}.
		\\
		
		(4) Use Theorem \ref{thm7} and the fact $(\eta(GMax(M)):_R M) \in GMax(R)$.
		\\
		
		(5) Since ${\eta}(\Upsilon_{M})=\{0\}\in \Upsilon_{M}$, by Theorem \ref{thm7}, $\Upsilon_{M}$ is an irreducible space.
		\\
		
		(6) Since $(0:_{R}M)=\{0\}\in GWSpec(R)$, so $\{0\}\in \Upsilon_{M}$
		(graded pseudo weakly prime submodule)\ by Theorem \ref{thm7}, $\Upsilon_{M}$ is an irreducible space.
	\end{corollary}
	
	\begin{lemma}\label{lem3}
		Let M be a G-graded pseudo weakly primeful R-module and let P be a graded radical ideal of R. Then $P=(PM:M)$ if and only if $Ann(M) \subseteq P$. In particular, $\omega M$ is a graded pseudo weakly primeful submodule of M for every $\omega \in \chi^R(Ann(M))$.
		\\
		\\
		\textbf{Proof.} Let P be a graded radical ideal of R and $Ann(M) \subseteq P = \bigcap\limits_i \omega_i$, where $\omega_i$ runs through $\chi^R(P)$. Since M is a G-graded pseudo weakly primeful R-module and $\omega_i \in \chi^R(Ann(M))$, then there exists a graded pseudo weakly primeful submodule $I_i$ of M such that $(I_i:M)=\omega_i$. Now, we have that $P \subseteq (PM:M)= ((\bigcap\limits_i \omega_i)M:M) \subseteq \bigcap\limits_i (\omega_iM:M) \subseteq (I_i:M)=\bigcap\limits_i \omega_i = P$. Therefore, $(PM:M)=P$. 
	\end{lemma}
	
	Consider the closed subset of a topological space $X$. If $X= Cl(\{x\})$, an element $x \in X$ is called a generic point of $X$. In Theorem \ref{thm4} (1), every element $J$ of $\Upsilon_M$ is a generic point of the irreducible closed subset $\chi(J)$, as we have seen. If the topological space is a $T_0$-space, a generic point of a closed subset $X$ of a topological space is unique (see Theorem \ref{thm4}). On $G$-graded $R$-modules, the next theorem is a nice application of the Zariski topology. Indeed, this shows that irreducible closed subsets of $\Upsilon_M$ and graded pseudo weakly prime submodules of the $G$-graded $R$-module $M$ has a relationship.
	
	\begin{theorem}\label{thm8}
		Let M be a G-graded R-module and $P \subseteq \Upsilon_M$.
		\\
		
		(1) Then P is an irreducible closed subset of $\Upsilon_M$ if and only if $P =\chi(I)$ for some $I \in \Upsilon_M$. As a result, there is a generic point in every irreducible closed subset of $\Upsilon_M$.
		\\
		
		(2) The correspondence $\chi (I) \rightarrow I$ is a bijection of the set of all irreducible components of $\Upsilon_M$ onto the set of all minimal elements of $\Upsilon_M$.
		\\
		
		(3) Let M be a G-graded pseudo weakly primeful R-module. Then the set of all irreducible components of $\Upsilon_M$ is of the form $H=\{\chi^M(LM):L$ is a minimal element of $\chi(Ann(M))\}$.
		\\
		\\
		\textbf{Proof.} (1) By Corollary \ref{coro2}, for any $I \in \Upsilon_M$, $P=\chi(I)$ is an irreducible closed subset of $\Upsilon_M$. Conversely, if P is an irreducible closed subset of $\Upsilon_M$, then for some $N\leq M$, $P = \chi (N)$ and $\eta(P)=\eta(\chi(N))= GPWrad(N) \in \Upsilon_M$ by Theorem \ref{thm7}. Thus $P=\chi(N)=\chi(GPWrad(N))$.
		\\
		
		(2) Let $\Upsilon_M$ have an irreducible component P. Each irreducible component of $\Upsilon_M$ is a maximum element of the collection $\{\chi (Q): Q \in \Upsilon_M\}$ by (1), then we have that $P=\chi(I)$ for some $I \in \Upsilon_M$. Obviously,$ I$ is a minimal element of $\Upsilon_M$, for if T is a graded pseudo weakly prime submodule of M with $T \subseteq I$, so $\chi(I)\subseteq \chi(T)$ and so that $I=T$. Now, Let I be a minimal element of $\Upsilon_M$ with $\chi(I)\subseteq\chi(Q)$ for some $Q\in \Upsilon_M$. Thus $Q = GPWrad(Q)=\eta(\chi(Q))\subseteq\eta(\chi(I))=GPWrad(I)=I$, therefore, $I = Q$. This means that the irreducible component of $\Upsilon_M$ is $\chi (I)$.
		\\
		
		(3) ($\Rightarrow$) Let $\Upsilon_M$ have an irreducible component P. By part (1), $P = \chi(I)$ for some $I \in \upsilon_M$. $(I: M)M$ is clearly a graded pseudo weakly prime submodule of M. Because $(I: M)M\subseteq I$, we have that $P=\chi(I)\subseteq \chi((I:M)M)$. Since P is an irreducible component, $\chi(I) = \chi((I:M)M)$, and so $I = (I:M)M$. We must prove that $i := (I:M)$ is a minimal element of $\chi^R(Ann(M))$. To see this, let $Y\in \chi^R(Ann(M))$ and $J\subseteq i$. Then $Y/Ann(M) \in GWSpec(R/Ann(M))$, and there exists an element $Q \in \Upsilon_M$ such that $(Q:M)=Y$ since M is a G-graded pseudo weakly primeful R-module. Hence $P=\chi^M(I) \subseteq \chi^M(Q)$. Thus $P = \chi^M(I)=\chi^M(Q)$ due to the maximality of $\chi^M(I)$ and we have $i=J$.  
		\\
		($\Leftarrow$) For $P \in H$, there exists a minimal element in $\chi^R(Ann(M))$ such that $P=\chi^M(LM)$. Since M is a G-graded pseudo weakly primeful R-module, LM is a graded pseudo weakly primeful submodule of M by Lemma \ref{lem3}. Thus P is an irreducible space by part (1). Assume that $P=\chi^M(LM)\subseteq \chi^M(Q)$, where Q is an element of $\Upsilon_M$. Since $LM \in \chi^M(Q)$ and L is a minimal element we have that $L=(LM:M)=(Q:M)$. Now, $P=\chi^M(LM)=\chi^M((Q:M)M)$ and $\chi^M(Q) \subseteq \chi^M((Q:M)M) $. Therefore, $P+ \chi^M(LM)=\chi^M(Q)$.
	\end{theorem}

	\begin{theorem}\label{thm9}
		Let M be a G-graded pseudo weakly primeful R-module, then the following statements are equivalent:
		\\
		
		(1) $\Upsilon_M$ is an irreducible space.
		\\
		
		(2) $GWSpes(R/Ann(M))$ is an irreducible space.
		\\
		
		(3) $\chi(Ann(M))$ is an irreducible space.
		\\
		
		(4) $Grad(Ann(M))$ is a graded weakly prime ideal of R.
		\\
		
		(5) $\Upsilon_M = \chi^M(PM)$ for some $P \in \chi(Ann(M))$.
		\\
		\\
		\textbf{Proof.} (1) $\Rightarrow$ (2) As we can observe from proof in Theorem \ref{thm3}, the natural map $\varphi$ is a continuous variable that is assumed to be surjective. Thus $Im(\varphi)=GWSpec(R/Ann(M))$ is also irreducible space.
		\\
		
		(2) $\Rightarrow$ (3) We define the mapping 
		\begin{center}
			$\phi : GWSpec(R/Ann(M)) \rightarrow GWSpec(R)$  \\
			$J/Ann(M) \ \rightarrow \ J$ 
		\end{center}
		which is a graded homomorphism. Therefore $\eta(Ann(M))$ is an irreducible space.
		\\
		
		(3) $\Rightarrow$ (4) By Theorem \ref{thm7}, we take $\eta(\chi(Ann(M)))= Grad(Ann(M))$ is a graded weakly prime ideal of R.
		\\
		
		(4) $\Rightarrow$ (5) By Lemma \ref{lem3}, $Grad(Ann(M))M$ is a graded pseudo weakly prime submodule of M. Now, let $I \in \Upsilon_M$, then $Grad(Ann(M)) \subseteq (I:M)$ and $Grad(Ann(M))M \subseteq I$. Thus $\Upsilon_M = \chi(Grad(Ann(M))M)$ where $Grad(Ann(M)) \in \chi(Ann(M))$.
		\\
		
		(5) $\Rightarrow$ (1) By Lemma \ref{lem3}, $P M$ is a graded pseudo weakly prime submodule of M and by Corollary \ref{coro2}, $\chi^M(PM)= \Upsilon_M$ is an irreducible space.
	\end{theorem}
	
	\begin{definition}\label{def3}
		Let R be a G-graded ring and T be a topological space of R.
		\\
		
		(1) T is said to be a graded spectral space if T is a homeomorphic to GSpec(R), with the zariski topology for R. 
		\\
		
		(2) T is said to be a graded weakly spectral space if T is a homeomorphic to GWSpec(R), with the zariski topology for R. 
	\end{definition}
	
	This work shows that the graded pseudo weakly prime spectrum of $G$-graded $R$-modules are a graded weakly spectral space under algebraic conditions. The topological space $T$ that satisfies the following conditions has been defined for graded weakly spectral spaces: (1) $T$ is a $T_0$-space, (2) $T$ is a quasi-compact, (3) Under finite intersection, the quasi-compact open subsets of $T$ are closed, (4) A generic point exists for any irreducible closed subset of $T$.
	
	\begin{remark}\label{rem3}
		A Noetherian topological space is a graded weakly spectral space if and only if it is a $T_0$-space and there is a generic point in every non-empty irreducible closed subspace.
	\end{remark}
	
	Recall that, if $M$ is a $G$-graded $R$-module, then $\Upsilon_M$ is a $T_0$-space (see Theorem \ref{thm4}) and every non-empty irreducible closed subset of $\Upsilon_M$ has a generic point (see Theorem \ref{thm8}).

	\begin{theorem}\label{thm10}
		Let M be a G-graded R-module. $\Upsilon_M$ is a graded weakly spectral space if R is a Noetherian graded ring and for every graded submodule P of M there exists a graded ideal I of R such that $\chi(P)=\chi(IM)$.
		\\
		\\
		\textbf{Proof.} We will prove that every open subset of $\Upsilon_M$ is quasi-compact. Let H be an open subset of $\Upsilon_M$ and $\{T_\delta\}_{\delta\in \Delta}$ be an open cover of H. Then there are a graded submodules P and $P_\delta$ of M such that $H = \Upsilon_M/\chi(P)$ and $T_\delta=\Upsilon_M/\chi(P_\delta)$ for some $\delta\in \Delta$ and $H \subseteq \bigcup\limits_{\delta\in\Delta} T_\delta = \Upsilon_M/\bigcap\limits_{\delta\in\Delta}\chi(P_\delta)$. Now, there exists a graded ideal $I_\delta$ of R such that $\chi(P_\delta)=\chi(I_\delta M)$. Thus, $H \subseteq \Upsilon/\chi((\sum\limits_{\delta\in\Delta} I_\delta) M)$. Since R is a Noetherian graded ring, there exist a finite subset $\Delta'$ of $\Delta$ such that $H \subseteq \bigcup\limits_{\delta\in \Delta'}T_\delta$. Therefore, $\Upsilon_M$ is a Noetherian topological space and whence a graded weakly spectral space.
	\end{theorem}

\bibliographystyle{amsplain}

\end{document}